\documentclass[12pt,draftcls,onecolumn]{IEEEtranv1}
\usepackage{cite}

\usepackage{header}
\newcommand{\pt}{u}
\newcommand{\y}{y}
\newcommand{\ev}{\Omega}
\usepackage{mathabx}

\newcommand{\vep}{\varepsilon}
\usepackage{amsmath, amssymb}
\begin{document}
%
\title{On the Characterization of Local Nash Equilibria in Continuous Games}

\newif\ifarxiv
\arxivtrue
\ifarxiv
\author{Lillian~J.~Ratliff,        Samuel~A.~Burden,and~S.~Shankar~Sastry,
\thanks{The authors are with the Department of Electrical Engineering and
Computer Sciences, University of California, Berkeley, Berkeley, CA, 94720 USA 
e-mail: $\{${\tt ratliffl, sburden, sastry}$\}${\tt@eecs.berkeley.edu}}
\thanks{This work is supported by NSF
CPS:Large:ActionWebs award number 0931843, TRUST (Team for Research in
Ubiquitous Secure Technology) which receives support from NSF (award
number CCF-0424422), FORCES (Foundations Of Resilient
CybEr-physical Systems) which receives support from NSF (award number
CNS-1239166).}
}
\else
\author{Lillian~J.~Ratliff,~\IEEEmembership{Student Member,~IEEE,}
        Samuel~A.~Burden,~\IEEEmembership{Member,~IEEE}
        and~S.~Shankar~Sastry,~\IEEEmembership{Fellow,~IEEE}
\thanks{The authors are with the Department of Electrical Engineering and
Computer Sciences, University of California, Berkeley, Berkeley, CA, 94720 USA 
e-mail: $\{${\tt ratliffl, sburden, sastry}$\}${\tt@eecs.berkeley.edu}}
\thanks{This work is supported by NSF
CPS:Large:ActionWebs award number 0931843, TRUST (Team for Research in
Ubiquitous Secure Technology) which receives support from NSF (award
number CCF-0424422), FORCES (Foundations Of Resilient
CybEr-physical Systems) which receives support from NSF (award number
CNS-1239166).}
}

\fi

\maketitle

\begin{abstract}
We present a unified framework for characterizing local Nash equilibria in
continuous games
on either infinite--dimensional or finite--dimensional non--convex strategy spaces. We
provide intrinsic necessary and sufficient first-- and second--order conditions
ensuring strategies constitute local Nash equilibria. We term points satisfying
the sufficient conditions \emph{differential Nash equilibria}. Further, we
provide a sufficient condition (non--degeneracy) guaranteeing differential Nash
equilibria are isolated and show that such equilibria are structurally stable.
We present tutorial examples to illustrate our results and highlight degeneracies that can arise in continuous games.\end{abstract}


\IEEEpeerreviewmaketitle

\section{Introduction}
\label{sec:intro}
\IEEEPARstart{M}{any} engineering systems are complex networks in which intelligent actors make
decisions regarding usage of shared, yet scarce, resources. 
Game theory provides established techniques for modeling competitive interactions
that have emerged as tools for analysis and synthesis of systems comprised of
dynamically--coupled decision--making agents possessing diverse and
oft--opposing interests~(see, e.g.~\cite{frihauf:2012aa,shamma:2005aa}). 
We focus on games with a finite number of agents where their strategy spaces are
continuous, either a finite--dimensional differentiable manifold or an
infinite--dimensional Banach manifold. 

Previous work on continuous games with convex strategy spaces and player costs
led to global characterization and computation of Nash
equilibria~\cite{Li:1987oq, Basar:1987kl,
Contreras:2004cr}. Adding constraints led to extensions of nonlinear
programming concepts, such as constraint qualification conditions, to games with
generalized Nash equilibria~\cite{dorsch:2013aa, Rosen:1965tg, facchinei:2007aa}.
Imposing a differentiable structure on the strategy spaces yielded other global
conditions ensuring existence and uniqueness of Nash equilibria and Pareto
optima~\cite{Thom:1974dk, Ekeland:1974pb, Smale:1975hc}. 
In contrast, we aim to analytically characterize and numerically compute
\emph{local} Nash equilibria in continuous games on non--convex strategy spaces.

Bounding the rationality of agents can result in \emph{myopic} behavior~\cite{Flam:1998kx}, meaning that agents seek strategies that that are optimal locally but not necessarily globally. 
Further, it is common in engineering applications for strategy spaces or player costs to be non--convex,
for example when an agent's configuration space is a constrained set or a differentiable manifold~\cite{Muhammad:2005zr, Klavins:2002ly}. 
These observations suggest that techniques for characterization and computation of local Nash equilibria have important practical applications.

Motivated by systems with myopic agents and non--convex strategy spaces, we seek
an intrinsic characterization for local Nash equilibria that is structurally stable and amenable to computation.  
By generalizing derivative--based conditions for local optimality in nonlinear programming~\cite{Bertsekas:1999fk} and optimal control~\cite{Polak:1997fk}, we provide necessary first-- and second--order conditions that local Nash equilibria must satisfy, and further develop a second--order sufficient condition ensuring player strategies constitute a local Nash equilibrium.
We term points satisfying this sufficient condition \emph{differential} Nash equilibria.
In contrast to a pure optimization problem, this second--order condition is
insufficient to guarantee a differential Nash equilibrium is isolated; in fact, 
games may possess a continuum of differential Nash equilibria. Hence, we introduce an additional second--order condition
ensuring a differential Nash equilibrium is isolated.


Verifying that a strategy constitutes a Nash equilibrium in non--trivial strategy spaces requires testing that a non--convex inequality constraint is satisfied on an open set, a task we regard as generally intractable. 
In contrast, our sufficient conditions for local Nash equilibria require only the evaluation of player costs and their derivatives at single points.
Further, our framework allows for numerical computations to be carried out when
players' strategy spaces and cost functions are non--convex. 
Hence, we provide tractable tools for characterization and computation of differential Nash equilibria in continuous games. 

We show that non--degenerate differential Nash equilibria are structurally
stable;
hence,  
measurement noise and modeling errors that give rise to a nearby game
 do not result in 
 drastically different equilibrium behavior---a property that is desirable in both the design of games
as well as inverse modeling of agent behavior in competitive environments. Further, structural stability
ensures that following the flow generated by the gradient of each player's cost converges locally to a stable, non--degenerate differential Nash equilibrium.
We remark that non--degenerate
differential Nash equilibria are generic in the
finite--dimensional case~\cite{ratliff:2014aa}.

The rest of the paper is organized as follows. 
In
section~\ref{sec:game} we present the game formulation in both the
finite--dimensional and infinite--dimensional case. We
follow with the characterization of local Nash equilibria in
Section~\ref{sec:characterization}. Throughout the paper we carry an example that provides
insight into the importance of the results and in Section~\ref{sec:example} we
return to the example in more detail. Finally, we conclude with discussion in
Section~\ref{sec:conclusion}. The necessary mathematical background and notation
is contained
in the Appendix.

\section{Game Formulation}
\label{sec:game}
The theory of games we consider concerns interaction between a finite number of rational agents generally having different interests and objectives. We refer to the rational agents as \emph{players}. Competition arises due to the fact that the players have opposing interests. 

Let us begin by considering a game in which we have $n$ selfish players with
competing interests. 
The strategy spaces are topological spaces $M_i$ for each $i\in\{1, \ldots,
n\}$. Note these can
be finite--dimensional smooth manifolds or infinite--dimensional Banach
manifolds. 
We denote the joint strategy space by $M=\prod_{i=1}^n M_i$.
The players are each interested in minimizing a cost function representing
their interests by choosing an element from their strategy space. We define
player $i$'s cost to be a twice--differentiable function $f_i\in C^2(M, \mb{R})$.
The following definition describes the equilibrium behavior we are interested
in:
\begin{definition}
  \label{def:SLNE}
  A strategy $(\pt_1,\ldots, \pt_n)\in M$ is a {\bfe{local Nash
  equilibrium}} if there exist open sets $W_i\subset M_i$ such
  that $\pt_i\in W_i$ and for each $i\in\{1,\ldots, n\}$ 
  \begin{equation}
    f_i(\pt_1, \ldots,\pt_i,\ldots, \pt_{n})\leq f_i(\pt_1, \ldots, \pt_i',\ldots, \pt_{n}),\ 
    \label{eq:ucondLN}
  \end{equation}
  for all  $\pt_i'\in
    W_i\bs\{\pt_i\}$. Further, if the above inequalities are strict, then we say
  $(\pt_1,\ldots, \pt_n)$ is a {\bfe{strict local Nash equilibrium}}.
  If $W_i=M_i$ for each $i$, then $(\pt_1,\ldots, \pt_n)$ is a {\bfe{global Nash
  equilibrium}}. 
\end{definition}
Simply put, the above definition says that no player can unilaterally deviate
from the Nash strategy and decrease their cost.

Before we move on to the characterization of local Nash equilibria, we
describe the types of games the results apply to and why they are important in
engineering applications.

Continuous games with finite--dimensional strategy spaces are described
by the player strategy spaces $M_1, \ldots, M_n$ and their cost functions
$(f_1,\ldots, f_n)$.
They
arise in a number of
engineering and economic applications, for instance, in modeling
one--shot decision making problems arising in 
transportation, communication and power networks~\cite{krichene:2014aa,
candogan:2010aa,park:2001aa}.
On the other hand, continuous games with infinite--dimensional strategy spaces, regarded as
open--loop differential games,
are used in engineering
applications in which there are agents coupled through dynamics. They
arise
in problems such as building energy
management~\cite{coogan:2013aa}, travel-time optimization in transportation
networks~\cite{bressan:2012aa}, and integration
of renewables into energy systems~\cite{zhu:2012aa}.

Open--loop differential games 
often come in the following form. 
Let $L_2[0,T]$ denote the space of square integrable functions from
$[0,T]\subset \mb{R}$ into $\mb{R}^m$. 
For an $n$--player game, strategy spaces are Banach manifolds, $M_i$ for $i\in\{1, \ldots,
n\}$, modeled on $L_2[0,T]$. For each $t\in [0,T]$, let $x(t)\in \mb{R}^n$ denote the state of the game. The state evolves according to the dynamics
\begin{equation}
  \dot{x}(t)=h(x(t),\uI{1}(t),\ldots, \uI{n}(t)) \ \ \forall \ t\in [0,T]
  \label{eq:statedynamics}
\end{equation}
where $\uI{i}\in M_i$ is player $i$'s strategy.  
We assume that $h(x, u_1,\ldots, u_n)$ is continuously differentiable,
globally Lipschitz continuous and all the derivatives in all its arguments are
globally Lipschitz
continuous.  We denote by $f_i(u_1, \ldots, u_n)=\hat{f}_i(x^{(x(0),\uI{1},\ldots, \uI{n})}(T))$ player
$i$'s cost function. The superscript notation on the state $x$ indicates the dependence of the state on the initial state and the strategies of the players. Each $\hat{f}_i$ is assumed twice continuously
differentiable so that 
each $f_i$ is 
$C^2$--Fr\'{e}chet--differentiable~\cite[Thm.~5.6.10]{Polak:1997fk}. We pose each player's optimization problem as
\begin{equation}
  \min_{\uI{i}}\hat{f}_i(x^{(x(0),\uI{1},\ldots, \uI{i}, \ldots, \uI{n})}(T)).
  \label{eq:optproblem}
\end{equation}
The costate for player $i$ evolves according to 
\begin{equation}
  \dot{p}_{i}(t)=-p_{i}(t)\pd{h}{x}(x(t), \uI{1}(t),\ldots,\uI{i}(t), \ldots,
  \uI{n}(t))
    \label{eq:costatedynamics}
\end{equation}
with final time condition 
\begin{equation}
  p_i(T)=D_xf_i(x^{(x(0),\uI{1}, \ldots, \uI{i},\ldots, \uI{n})}(T)).
  \label{eq:finaltime}
\end{equation}
The derivative of the $i$--th player's cost function is given by
\begin{equation}
  (D_{i}f_{i})(t)=p_i(t)\pd{h}{\uI{i}}(x(t), \uI{1}, \ldots, \uI{i}(t),\ldots,
  \uI{n}(t)).
  \label{eq:Dufi}
\end{equation}

Before we dive into the details, let us consider a simple example that exhibits
very interesting behavior.
\begin{example}[Betty--Sue]
  \label{ex:bettysue}
 Consider a two player game between Betty and Sue.  Let Betty's strategy
 space be $M_1=\mb{R}$ and her cost function
 $f_1(u_1,u_2)=\frac{u_1^2}{2}-u_1u_2$.
 Similarly, let Sue's strategy space be  $M_2=\mb{R}$ and her cost function
 $f_2(u_1,u_2)=\frac{u_2^2}{2}-u_1u_2$. 
 This game can be thought of as an abstraction of two agents in a building
 occupying adjoining rooms. The first term in each of their
 costs represents an energy cost and the second term is a cost from
 thermodynamic coupling. The agents try to maintain the temperature at a desired
 set--point in thermodynamic equilibrium. 

 \newcommand{\vv}{p}
 \newcommand{\w}{q}
Definition~\ref{def:SLNE} specifies that a point $(\vv,\w)$ is a Nash equilibrium if no player can unilaterally
deviate and decrease their cost, i.e. $f_1(\vv,\w)<f_1(u_1,\w)$
for all $u_1\in\mb{R}$
and $f_2(\vv,\w)<f_2(\vv,u_2)$ for all $u_2\in\mb{R}$.
  
Fix Sue's strategy $u_2=\w$, and calculate
  \begin{equation}
    D_1 f_1 = \pd{f_1}{u_1}=u_1-\w
    \label{eq:expdf1}
  \end{equation}
  Then, Betty's optimal response to Sue playing $u_2=\w$ is $u_1=\w$. Similarly, if
  we fix $u_1=\vv$, then Sue's optimal response to Betty playing $u_1=\vv$ is
  $u_2=\vv$. For all $u_1\in \mb{R}\bs\{\w\}$
  \begin{equation}
    -\frac{\w^2}{2}<\frac{u_1^2}{2}-u_1\w
    \label{eq:ineq1}
  \end{equation}
  so that $f_1(\w,\w)<f_1(u_1,\w)$ for all $u_1\in \mb{R}\bs\{\w\}$. Again,
  similarly, for all $u_2\in \mb{R}\bs\{\vv\}$
  \begin{equation}
    -\frac{\vv^2}{2}<\frac{u_2^2}{2}-u_2\vv
    \label{eq:ineq2}
  \end{equation}
  so that $f_2(\vv,\vv)<f_2(\vv,u_2)$ for all $u_2\in \mb{R}\bs\{\vv\}$. Hence, all the
  points on the line $u_1=u_2$ in $M_1\times M_2=\mb{R}^2$ are strict local Nash
  equilibria---in fact, they are strict global Nash equilibria.  \hfill\ensuremath{\blacksquare}
\end{example}

As the above example shows, continuous games can exhibit a continuum of equilibria.
Throughout the text we will return to this example.

\section{Characterization of Local Nash Equilibria}
\label{sec:characterization}
In this section, we characterize local Nash equilibria by
paralleling results in nonlinear programming and optimal control that provide
first-- and second--order necessary and sufficient conditions for local optima.

The following definition of a differential game form is due to Stein~\cite{Stein:2010fe}. 

\begin{definition}
  A \bfe{differential game form} is a differential $1$--form $\omega:
  M_1\times\cdots\times M_n\rar T^\ast(M_1\times\cdots\times M_n)$ defined by
  \begin{equation}
    \omega=\sum_{i=1}^n\psi_{M_i}\circ \dm f_i. 
    \label{eq:gf}
  \end{equation}
  where $\psi_{M_i}$ are the natural bundle maps defined in
  \eqref{eq:bundlemaps} that annihilate those components of the covector $\dm f_i$
  not corresponding to $M_i$.
  \label{def:gameform}
 \end{definition}
 \begin{remark}
   If each $M_i$ is a finite--dimensional manifold of dimension $m_i$, then the differential
   game form has the following coordinate representation:
   \begin{equation}
     \omega_{\vphi}=\sum_{i=1}^n\sum_{j=1}^{m_i}\pd{(f_i\circ\vphi^{-1})}{y^j_i}dy^{j}_i 
    \label{eq:gfcoord}
  \end{equation}
  where $(U, \vphi)$ is a product chart on $M$ at $\pt=(\pt_1, \ldots,
  \pt_n)$  with local coordinates $(\ucn{1}{1}, \ldots,$ $\ucn{m_1}{1},\ldots,
  \ucn{1}{n},
  \ldots, \ucn{m_n}{n})$ and where $U=\prod_{i=1}^nU_i$
  and $\vphi=\bigtimes_{i=1}^n\vphi_i$.  In addition, $f_i\circ\vphi^{-1}$ is the
  coordinate representation of $f_i$ for $i\in\{1,\ldots, n\}$.
 In particular, $\vphi_i(\pt_i)=(\ucn{1}{i},
  \ldots, \ucn{m_i}{i})$ where each $\ucn{j}{i}:U_i\rar \mb{R}$ is a
  coordinate function so that $d\ucn{j}{i}$ is its derivative.\hfill\ensuremath{\blacksquare}
\end{remark}
The differential game form captures a differential view of
the strategic interaction between the players. Indeed, $\omega$ indicates the
direction in which the players can change their strategies to decrease their
individual cost functions most rapidly. Note that each player's cost function
depends on its own choice variable as well as all the other players' choice
variables. However, each player can only affect their payoff by adjusting
their own strategy.

\begin{definition}
  \label{def:DNE}
  A strategy $u=(\pt_1,\ldots,\pt_n)\in M_1\times \cdots\times M_n$ is a \bfe{differential Nash equilibrium}
  if $\omega(\pt)=0$  and  $\D^2_{ii}f_i(\pt)$ is
  positive--definite for each $i\in\{1, \ldots,  n\}$.
\end{definition}
The second--order conditions used to define differential Nash equilibria are
motivated by results in nonlinear programming that use first-- and second--order
conditions to assess whether a critical point is a local optima \cite{Polak:1997fk}, \cite{Bertsekas:1999fk}.

The following proposition provides first-- and second--order necessary conditions for local Nash equilibria.
We remark that these conditions are reminiscent of those seen in nonlinear programming for optimality of critical points.
\begin{prop}
  If $u=(\pt_1,\ldots, \pt_n)$ is a local Nash equilibrium, then 
  $\omega(\pt)=0$ and
  $\D_{ii}^2f_i(\pt)$ is positive semi-definite for each
  $i\in\{1,\ldots, n\}$.
\end{prop}
\begin{IEEEproof}
  Suppose that $u=(\pt_1,\ldots, \pt_n)\in M$ is a local Nash equilibrium. Then, 
  \begin{equation}
    f_i(\pt)\leq f_i(\pt_1,\ldots,\pt_i',\ldots, \pt_{n}),\ \ \forall \ \pt_i'\in
    W_i\bs\{\pt_i\}
    \label{eq:ucondLN-p}
  \end{equation}
  for open $W_i\subset M_i$, $i\in\{1,\ldots, n\}$.
  Suppose that we have a product chart $(U, \vphi)$, where $U=\prod_{i=1}^nU_i$ and $\vphi=\bigtimes_{i=1}^n\vphi_i$, such that
  $\pt\in U$. 
  
  Let $\vphi_i(\pt_i)=v_i$ for each $i$. 
  Then, since $\vphi$ is continuous, for each $i\in\{1,\ldots, n\}$, we have that for all $v_i'\in
  \vphi_i(W_i\cap U_i)\bs\{\vphi_i(\pt_i)\}$,
  \begin{equation}
    f_i\circ \vphi^{-1}(v_1,\ldots,v_i, \ldots, v_{n})\leq f_i\circ
    \vphi^{-1}(v_1, \ldots, v_i',\ldots, v_{n}).
    \label{eq:local-cond1}
  \end{equation}
  Now, we apply Proposition 1.1.1 from
  \cite{Bertsekas:1999fk},
  if $M_i$ is finite--dimensional, or Theorem
  4.2.3(1) and Theorem 4.2.4(a) from \cite{Polak:1997fk}, if $M_i$ is
  infinite--dimensional, to $f_i\circ
  \vphi^{-1}$. We conclude that for each $i\in\{1,\ldots,n\}$,
 $\D_i(f_i\circ
  \vphi^{-1})(v_1,\ldots,v_{n})=0$
   and for all $\nu \in \vphi_i(U_i\cap W_i)$,
   \begin{equation}
     D_{ii}^2(f_i\circ\vphi^{-1})(v_1,\ldots, v_{n})(\nu, \nu)\geq
     \alpha \|\nu\|^2, 
     \label{eq:hesscond}
   \end{equation}
 i.e. it is a
  positive semi--definite bilinear form on $\vphi_i(U_i\cap W_i)$.
  
  Invariance of the stationarity of critical points and the index of the Hessian with
  respect to coordinate change gives us $\omega(\pt)=0$ and
  $\D^2_{ii}f_i(\pt)$ is a positive semi--definite for each $i\in\{1,\ldots,
  n\}$. 
\end{IEEEproof}

We now show that the conditions defining a differential Nash equilibrium are sufficient to guarantee a strict local Nash equilibrium. 
\begin{thm}
  A differential Nash equilibrium is a strict local Nash equilibrium.
  \label{prop:DNE-SLNE}
\end{thm}
\begin{IEEEproof}
  Suppose that $u=(\pt_1,\ldots, \pt_n)\in M$ is a differential Nash equilibrium.
  Then, by the definition of differential Nash equilibrium, $\omega(\pt)=0$ and
  $\D^2_{ii}f_i(\pt)$ is positive definite for each
  $i\in\{1, \ldots, n\}$. The second-derivative conditions imply that  $\D_{ii}^2(f_i\circ \vphi^{-1})(v_1, \ldots, v_n)$ is a
  positive--definite bilinear form where $v_i=\vphi_i(u_i)$ for any
  coordinate chart $(U, \vphi)$, with $\vphi=\bigtimes_{i}\vphi_i$,
  $U=\prod_i U_i$, and $\pt_i\in U_i$ for each $i\in \{1, \ldots, n\}$.
  

Using the isomorphism introduced in the appendix in~\eqref{eq:canon}, 
$\omega(\pt)=0$ implies that for each $i\in\{1,\ldots, n\}$,
    $\D_i(f_i\circ \vphi^{-1})(v_1,\ldots, v_n)=0.$
Let $E_i$ be the model space, i.e. the underlying
Banach space, in either the finite--dimensional or infinite--dimensional case.
Applying either Proposition~1.1.3 from \cite{Bertsekas:1999fk} or Theorem 4.2.6
(a) from \cite{Polak:1997fk} to to each $f_i\circ \vphi^{-1}$ with
$(\vphi_{1}(\pt_{1}), \ldots, \vphi_{i-1}(\pt_{i-1}), \vphi_{i+1}(\pt_{i+1}),
\ldots, \vphi_n(u_n))$ fixed yields
a neighborhood
$W_i\subset E_i$ such that for all $v'\in W_i$, 
\begin{equation}
  f_i\circ \vphi^{-1}(v_1,\ldots,v_i, \ldots, v_n)<f_i\circ
  \vphi^{-1}(v_1,\ldots,v', \ldots,v_{n}).
  \label{eq:nasheqi}
\end{equation}
Since $\vphi$ is continuous, there exists a neighborhood $V_i\subset M_i$
of $\pt_i$ such that for $V_i=\vphi_i^{-1}(W_i)$ and all $\pt_i'\in V_i\backslash \{\pt_i\}$, 
\begin{equation}
  f_i(\pt_1,\ldots,\pt_i,\ldots, \pt_{n})<f_i(\pt_1,\ldots,\pt_i',\ldots,
  \pt_{n}).
  \label{eq:fnashM}
\end{equation}

%

Therefore, differential Nash equilibria are strict local Nash equilibria.
Due to
the fact that both $\omega(u)=0$ 
and definiteness of the
Hessian are coordinate invariant, this is independent of choice of coordinate
chart.
  \end{IEEEproof}

 We remark that the conditions for differential Nash equilibria are not
  sufficient to guarantee that an equilibrium is isolated. 
    \addtocounter{example}{-1}
  \begin{example}[Betty--Sue: Continuum of Differential Nash]
    Returning to the Betty--Sue {\debacle}, we can check that at all the points such
    that
    $u_1=u_2$, $\omega(u_1,u_2)=0$ and $D_{ii}^2f_i(u_1, u_2)=1>0$ for each
    $i\in\{1,2\}$. Hence, there is a continuum of differential Nash equilibria
    in this game.\hfill\ensuremath{\blacksquare}
  \end{example}
We propose a sufficient condition to guarantee that differential Nash equilibria
are isolated.
We do so by combining ideas introduced by Rosen for convex games with concepts
from Morse theory, in particular second--order conditions on non--degenerate critical points of real-valued functions on manifolds. 

At a differential Nash equilibrium $\pt=(\pt_1,
  \ldots, \pt_n)$, consider the derivative of the differential game form
\begin{align}
  d\omega= \sum_{i=1}^nd(\psi_{M_i} \circ df_i).
  \label{eq:Domega}
\end{align}
  Intrinsically, this derivative is a tensor field $d\omega\in T^0_2(M)$; at a point $\pt\in M$ where $\omega(\pt) = 0$ it is a bilinear form constructed
  from the uniquely determined continuous, symmetric, bilinear forms
  $\set{d^2f_i(\pt)}_{i=1}^n$.


\begin{thm}\label{thm:Dw}
  If $\pt=(\pt_1,
  \ldots, \pt_n)$ is a differential Nash equilibrium and
  $d\omega(\pt)$ is non--degenerate, then $\pt$ is an isolated strict local Nash equilibrium.   
\end{thm}
\begin{IEEEproof}
 Since $\pt$ is a differential Nash equilibrium, Theorem~\ref{prop:DNE-SLNE} gives us that it is a strict local Nash equilibrium. The following argument shows that it is isolated.
 Choose a coordinate chart $(U, \vphi)$ with $\vphi=\bigtimes_{i=1}^n\vphi_i$
 and $U=\prod_{i=1}^n U_i$. Let $E$ denote the underlying model space of the
 manifold $M_1\times \cdots \times M_n$. Define the map $g: E \rar E$ by
 \begin{equation}
   g(\vphi(\pt))=\sum_{i=1}^n D_i(f_i\circ \vphi^{-1})(\vphi(\pt))
   \label{eq:gmap}
 \end{equation}
Note that $g$ is the coordinate representation of the differential game
form $\omega$. Zeros of the function $g$ define critical points of the
game and its derivative at critical points is $d\omega$. Since $\pt$ is a
differential Nash equilibrium, $\omega(\pt)=0$. Further, since
$d\omega(\pt)$ is non--degenerate---the map $A(v)(w)=d\omega(\pt)(v,w)$
is a linear isomorphism---we can apply the Inverse Function   
 Theorem~\cite[Thm.~2.5.2]{abraham:1988aa} to get that $g$ is a local
 diffeomorphism at $\pt$, i.e. there exists an open neighborhood $V$ of
 $\pt$ such that the restriction of $g$ to $V$ establishes a diffeomorphism
 between $V$ and an open subset of $E$. Thus, only $\vphi(\pt)$ could be
 mapped to zero near $\vphi(\pt)$. Non--degeneracy of $d\omega(\pt)$ is
 invariant with respect to choice of coordinates. Therefore, $\pt$ is isolated.
\end{IEEEproof}
 
\begin{definition}
  Differential Nash equilibira $\pt=(\pt_1, \ldots, \pt_n)$ such that
  $d\omega(\pt)$ is non--degenerate are termed \emph{non--degenerate
  differential Nash equilibria}. 
\end{definition}

  \addtocounter{example}{-1}
\begin{example}[Betty--Sue: Degeneracy and Breaking Symmetry]
  Return again to the Betty--Sue example in which we showed
  that there is a continuum of Nash equilibria; in fact, all the points on the
  line $u_1=u_2$ are differential Nash equilibria and at each of these points we
  have
%
  \begin{equation}
    d\omega(u_1,u_2)=\bmat{\ \ 1 & -1\\ -1& \ \ 1}
    \label{eq:domegex}
  \end{equation}
  so that $\det(d\omega(u_1,u_2))=0$. Hence, all of the equilibria are
  \emph{degenerate}. 
  By breaking the symmetry in the game, we can make $(0,0)$ a
  non--degenerate differential Nash equilibrium; i.e. we can remove all but one
  of the equilibria. Indeed, let Betty's cost be
  given by $\tilde{f}_1(u_1, u_2)=\frac{u_1^2}{2}-au_1u_2$ and let Sue's cost
  remain unchanged. Then the local representation of the derivative of the differential
  game form $\td{\omega}$ of the game $(\td{f}_1, f_2)$ is 
  \begin{equation}
    d\td{\omega}(u_1,u_2)=\bmat{\ \ 1 & -a\\ -1& \ \ 1}
    \label{eq:domegex_}
  \end{equation}
Thus for any value of $a\neq 1$, $(0,0)$ is a non--degenerate differential Nash
equilibrium. This shows that small modeling errors can remove degenerate
differential Nash equilibria.
\hfill\ensuremath{\blacksquare}

\end{example}

In a neighborhood of a non--degenerate differential Nash equilibrium there are no other Nash equilibria. This
property is desirable particularly in applications where a central planner is
designing incentives to induce a socially optimal or otherwise desirable equilibrium that optimizes the central
planner's cost; if the desired equilibrium
resides on a continuum of equilibria, then due to measurement noise or myopic
play, agents may be induced to play a nearby equilibrium that is 
suboptimal for the central planner. In Section~\ref{sec:example}, we extend
Example~\ref{ex:bettysue} by introducing a central planner.
But first, we show that non--degenerate differential Nash equilibria are structurally stable.

\section{Structural Stability}
\label{sec:structuralstability}
Examples demonstrate that global Nash equilibria may fail to persist under
arbitrarily small changes in player costs~\cite{Ekeland:1974pb}. A natural question arises: do local Nash equilibria persist under perturbations? Applying structural stability analysis from dynamical systems theory, we answer this question affirmatively for non–degenerate differential Nash equilibria subject to smooth perturbations in player costs.

Let $M=M_1\times\cdots\times M_n$ and $f_1,\ldots,f_n:M\into\R$ be $C^2$ player cost functions,
$\omega:M\into T^\ast M$ the associated differential game form \eqref{eq:gf},
and suppose $\pt\in M$ is a non--degenerate differential Nash equilibrium,
i.e. $\omega(\pt) = 0$ and $d\omega(\pt)$ is non--degenerate.
We show that 
for all $\td{f}_i\in
C^\infty(M,\R)$ sufficiently close to $f_i$ there exists a unique
non--degenerate differential Nash equilibrium $\td{\pt}\in M$ for
$(\td{f}_1,\ldots,\td{f}_n)$ near $\pt$.

\begin{prop}[Parameterized Structural Stability] \label{thm:stable_param}
Non--degenerate differential Nash equilibria are parametrically structurally stable: 
given $f_1,\ldots, f_n\in C^2(M,\R)$, $\zeta_1,\ldots, \zeta_n\in C^2(M,\R)$,
and a non--degenerate differential Nash equilibrium $\pt\in M$ for
$(f_1,\ldots, f_n)$,
there exist neighborhoods $U\subset\R$ of $0$ and $W\subset M$ of $\pt$
such that for all $s\in U$ there exists a unique non--degenerate differential
Nash equilibrium $\td{\pt}(s)\in W$ for $\paren{f_1 + s\zeta_1,\ldots, f_n
+ s\zeta_n}$.
\end{prop}

\begin{IEEEproof}
Define $\td{f}_j: M_1\times \cdots \times M_n\times\R\into\R$ by 
\[ \td{f}_j(\pt,s) = f_j(\pt) + s\zeta_j(\pt) \]
and $\td{\omega}: M_1\times \cdots \times M_n\times\R\into T^\ast (M_1\times
\cdots \times M_n)$ by
\[
\td{\omega}(\pt,s) = \sum_{i=1}^n \td{\psi}_{M_i}\circ d\td{f}_i(\pt,s)
\]
for all $s\in\R$ and $\pt\in M_1\times \cdots \times M_n$ and where $\td{\psi}_{M_i}:T^\ast(M_1\times
\cdots\times M_n\times \mb{R})\rar T^\ast(M_1\times
\cdots\times M_n\times \mb{R})$.
Observe that $\D_{1}\td{\omega}( (u_1, \ldots, u_n),0)$ is invertible since $\pt$ is a
non--degenerate differential Nash equilibrium for $(f_1,\ldots,f_n)$.
Therefore by the Implicit Function Theorem~{\cite[Prop. 3.3.13~(iii)]{abraham:1988aa}}, there exist neighborhoods $V\subset\R$ of $0$ and
$W\subset M$ of $\pt$ and a smooth function $\sigma\in C^\infty(V,W)$ such that 
\[ \forall s\in V,\pt\in W : \td{\omega}(\pt,s) = 0 \iff \pt = \sigma(s). \]

Furthermore, since $\td{\omega}$ is continuously differentiable, there exists a
neighborhood $U\subset V$ of $0$ such that $d\td{\omega}(\sigma(s),s)$ is invertible for all $s\in U$.
We conclude for all $s\in U$ that $\sigma(s)\in M$ is the unique Nash
equilibrium for $\paren{\paren{f_1 + s\zeta_1}|_W,\ldots,\paren{f_n + s\zeta_n}|_W}$, and furthermore that $\sigma(s)$ is a non--degenerate differential Nash equilibrium.
\end{IEEEproof}

We remark that the preceding analysis extends directly to any finitely--parameterized perturbation.
For an arbitrary perturbation, we have the following.

\begin{thm}[Structural Stability] \label{thm:stable}
Non--degenerate differential Nash equilibria are structurally stable: 
let $\pt\in M$ be a non--degenerate differential Nash equilibrium
for $(f_1,\ldots, f_n)\in C^2(M,\R^n)$.
Then there exist neighborhoods $U\subset C^2(M,\R^n)$ of $(f_1,\ldots,
f_n)$ and $W\subset M$ of $\pt$ and a $C^2$ Fr\'{e}chet--differentiable function
$\sigma\in C^2(U,W)$ such that 
for all $(\td{f}_1,\ldots,\td{f}_n)\in U$ the point $\sigma(\td{f}_1,\ldots,
\td{f}_n)$ is the unique non--degenerate differential Nash equilibrium for
$(\td{f}_1,\ldots, \td{f}_2)$ in $W$.
\end{thm}

\begin{IEEEproof}
Consider the operator
$\ev\in C^1(C^1(M,\R^n)\times M, \R^n)$
 defined by 
\begin{align}
  \ev((\td{f}_1,&\ldots, \td{f}_n),(\pt_1, \ldots, \pt_n)) =\sum_{i=1}^n \psi_{M_i}\circ
  d\td{f}_i(\pt_1, \ldots, \pt_n).
  \label{eq:eval}
\end{align}
Note that the right--hand side is the differential game form $\td{\omega}(\pt_1,
\ldots, \pt_n)$ for the game $(\td{f}_1, \ldots, \td{f}_n)$. Suppose that $u=(u_1, \ldots, u_n)$ is a non-degenerate
differential Nash equilibrium. A straightforward application of
Proposition~2.4.20~\cite{abraham:1988aa} implies that the operator $\ev$ is
$C^1$
Fr\'{e}chet--differentiable. In addition,
\begin{equation}
  \D_2\ev((f_1,\ldots, f_n),(\pt_1, \ldots, \pt_n))=d\omega(\pt_1,
  \ldots, \pt_n).
  \label{eq:dev}
\end{equation}
Since $d\omega(\pt)$ is an isomorphism by assumption, we can apply the Implicit
Function Theorem \cite[Prop.~3.3.13~(iii)]{abraham:1988aa} to $\ev$ to get an open neighborhood $W\subset M$ of $\pt$
and $V\subset C^2(M, \mb{R}^n)$ of $(f_1, \ldots, f_n)$ and a smooth
function $\sigma\in C^2(V, W)$ such that
$$\forall \tilde{f}\in V, \ v\in W:\  \ev(\tilde{f}, v)=0 \Longleftrightarrow \
v=\sigma(\tilde{f})$$
where $\tilde{f}=(\tilde{f}_1, \ldots, \tilde{f}_n)$. Furthermore, since $\ev$
is continuously differentiable, there exists a neighborhood $U\subset V$ of
$(f_1, \ldots, f_n)$ such that $d\ev(\tilde{f}, \sigma(\tilde{f}))$ is
invertible for all $\tilde{f}\in U$. Thus, for all $\tilde{f}\in U$,
$\sigma(\tilde{f})\in M$ is the unique non--degenerate differential Nash
equilibrium.
\end{IEEEproof}

Let us return to
Example~\ref{ex:bettysue} and examine what can happen in the degenerate case.
 \addtocounter{example}{-1}

\begin{example}[Betty--Sue: Structural Instability]
 Let us recall again the Betty--Sue {\debacle} in which we have a game admitting a continuum of differential Nash equilibria. We can show that an arbitrarily small perturbation will make \emph{all} the equilibria disappear. 
Indeed, let $\vep\neq 0$ be arbitrarily small and consider Betty's perturbed cost function
\begin{equation}
  \tilde{f}_1(u_1,u_2)=\frac{u_1^2}{2}-u_1u_2+\vep u_1.
  \label{eq:perturbedUrbain}
\end{equation}
Let Sue's cost function remain unchanged. Then, \emph{all} Nash equilibria
disappear. Indeed, a necessary condition that a Nash equilibrium $(u_1,u_2)\in
M_1\times M_2$ must satisfy is $\omega(u_1,u_2)=0$ thereby implying
$\D_1\tilde{f_1}(u_1,u_2)=u_1-u_2+\vep=0$ and $\D_2f_2(u_1,u_2)=u_2-u_1=0$. This can
only happen for $\vep=0$. Hence, \emph{any} perturbation $\vep u_1$ with
$\vep\neq 0$ will remove all the Nash equilibria. \hfill\ensuremath{\blacksquare}

\end{example}

On the other hand, equilibria
that are stable---thereby attracting using decoupled myopic approximate
best–-response---persist under small perturbations~\cite{ratliff:2013aa}.

\begin{example}[Convergence of Gradient Play]
  \label{ex:gradplay}
We adopt a dynamical systems perspective of a two--player game over the strategy space $U_1\times U_2$ with player costs $f_1,f_2:U_1\times U_2\into\R$.
Specifically, we consider the continuous--time dynamical system generated by the negative of the player's individual gradients:
\eqnn{\label{eq:ct}
\mat{c}{\dot{u}_1 \\ \dot{u}_2} = \mat{c}{- D_1 f_1(u_1,u_2) \\ - D_2 f_2(u_1,u_2)} = -\omega(u).
}
If $(\mu_1,\mu_2)\in U_1\times U_2$ is a differential Nash equilibrium, then $\omega(\mu_1,\mu_2) = 0$. 
These dynamics are \emph{uncoupled} in the sense the dynamics  $\dot{u}_i$ for each player do not depend on the cost function of the other player. 
It is known that such uncoupled dynamics need not converge to local Nash equilibria~\cite{hart:2003aa}.
However, the subset of
non--degenerate differential Nash equilibria where the
spectrum of  $d\omega$ is strictly in the right--half plane (in the
finite--dimensional case, this corresponds to 
all eigenvalues of $d\omega$
having strictly positive real parts) are exponentially stable stationary points of
\eqref{eq:ct}~\cite[Prop.~4]{ratliff:2013aa},~\cite[Thm.~4.3.4]{abraham:1988aa}.
Theorem~\ref{thm:stable} shows that convergence of uncoupled gradient play to
such \emph{stable} non--degenerate differential Nash equilibria persists under
small smooth perturbations to player costs. \hfill\ensuremath{\blacksquare}

\end{example}

We remark that in the finite--dimensional case we can show that non--degenerate
differential Nash equilibria are generic among local Nash
equilibria~\cite{ratliff:2014aa}. Genericity implies that local Nash equilibria
in an open--dense set of continuous games are non--degenerate differential Nash
equilibria. 
Furthermore, structural stability implies that these equilibria
persist under smooth perturbations to player costs.
As a consequence,
small modeling errors or environmental disturbances generally do not result in
games with drastically different equilibrium behavior.

\section{Inducing a Nash Equilibrium}
\label{sec:example}


The problem of inducing Nash equilibria through incentive mechanisms appears in
engineering applications including energy management~\cite{coogan:2013aa} and
network security~\cite{Ratliff:2012uq, zhu:2012ac}. The
central planner aims to shift the Nash equilibrium of the agents' game to one
that is desirable from its perspective. 
Thus the central planner 
optimizes its cost subject to
constraints 
given by the inequalities that define a Nash equilibrium.
This
requires verification of non--convex conditions on an open set---a generally intractable task. A
natural solution is to replace these inequalities with first-- and second--order
 sufficient conditions on each agent's optimization problem.
As the Betty--Sue {\debacle} shows (Example~\ref{ex:bettysue}), these necessary conditions are
not enough to guarantee the desired Nash is isolated; the additional constraint
that $d\omega$ be non--degenerate
must be enforced.


 \addtocounter{example}{-2}

 \begin{example}[Betty--Sue: Inducing Nash]
   Consider a central planner who desires to optimize the cost of deviating from
   the temperature $\tau$:
\begin{equation}
 f_p(u_1,u_2)=(u_1-\tau)^2+(u_2-\tau)^2.
  \label{eq:spcpst}
\end{equation}
The central planner wants to induce the agents to play $(u_1,u_2)=(\tau,\tau)$ 
by selecting $a\in \mb{R}$ and augmenting Betty's and Sue's costs:
$$\td{f}_1^a(u_1,u_2)=f_1(u_2,u_2)+\frac{a}{2}(u_1-\tau)^2$$ 
$$\td{f}_2^a(u_1,u_2)=f_2(u_1,u_2)+\frac{a}{2}(u_2-\tau)^2.$$
The differential game form of the augmented game $(\td{f}_1^a, \td{f}_2^a)$ is
$$\td{\omega}(u_1,u_2)=(u_1-u_2+a(u_1-\tau))du_1+(u_2-u_1+a(u_2-\tau))du_2$$
and the second--order differential game form is 
$$d\td{\omega}(u_1,u_2)=\bmat{1+a & -1\\ -1 & 1+a}.$$
For any $a\in(-1, \infty)$, $(\tau,\tau)$ is a differential Nash equilibrium of
$(\td{f}_1^a, \td{f}_2^a)$ since $\td{\omega}(\tau,\tau)=0$ and
$d_{ii}^2\td{f}_i^a(\tau,\tau)>0$. For
any $ a\in(-1,0]$, the game $(\td{f}_1^a,
\td{f}_2^a)$ undesirable behavior. Indeed, recall Example~\ref{ex:gradplay} in which we
consider the gradient dynamics for a two player game. For values of $a\in (-1,
0)$, $d\td{\omega}$ is indefinite so that the equilibrium of the gradient system is a
saddle point. Hence, if agents perform gradient play and happen to initialize
on the unstable manifold, then they will not converge to any equilibrium. Further, while $a=0$
seems like a natural choice since it means not augmenting the players costs at
all, it in fact gives rise to a continuum of equilibria. 
However, for $a>0$,
$d\td{\omega}$ is positive definite so that, as Example~\ref{ex:gradplay} points out,
the gradient dynamics will converge and the value of $a$ determines the
contraction rate.\hfill \ensuremath{\blacksquare}
\end{example}


This example indicates how
undesirable behavior can arise when the operator $d\omega$ is degenerate. Further, if the goal is to induce a particular Nash
equilibrium amongst competitive agents, then it is not enough to consider only
 necessary and sufficient conditions for Nash equilibria; inducing stable non--degenerate differential Nash equilibria
leads to desirable and structurally stable behavior.

\section{Discussion}
\label{sec:conclusion}
By paralleling results in non--linear programming and optimal control, we
developed first-- and second--order necessary and sufficient conditions that
characterize local Nash equilibria in continuous games on both finite-- and
infinite--dimensional strategy spaces. We further provided a second--order
sufficient condition guaranteeing differential Nash equilibria are
non--degenerate and, hence, isolated. 
We showed that non--degenerate differential
Nash equilibria are structurally stable and thus small modeling errors or
environmental disturbances generally will not result in games with drastically different equilibrium behavior.
Further, as a result of structural stability, our characterization of non--degenerate differential Nash equilibria is
amenable to computation. 
We illustrate through an example that such a characterization has
value for the design of incentives to induce a desired equilibria. By enforcing
not only non--degeneracy but also stability of a differential Nash equilibrium, the
central planner can ensure that the desired equilibrium is isolated and that
gradient play will  converge locally.



%

%
%
%
%
%
%
\appendix[Mathematical Prelimiaries]
\label{sec:prelims}
This appendix contains the standard mathematical objects used throughout this
paper (see~\cite{Lee:2012fy, abraham:1988aa} for a more detailed introduction).    

Suppose that $M$ is second--countable and a Hausdorff topological space. Then a \emph{chart} on $M$ is a 
homeomorphism $\vphi$ from an open subset $U$ of $M$ to an open subset of a Banach space.
We sometimes denote a chart by the
pair $(U,
\vphi)$.
Two charts $(U_1, \vphi_1)$ and $(U_2,
\vphi_2)$ are \emph{$C^r$--compatible} if and only if the composition $\vphi_2\circ
\vphi_1^{-1}:\vphi_1(U_1\cap U_2)\rar \vphi_2(U_1\cap U_2)$ is a 
$C^r$--diffeomorphism. A \emph{$C^r$--atlas} on $M$ is a collection of charts
$\{(U_\alpha, \vphi_\alpha)\}_{\alpha\in\mc{A}}$ any two of which are $C^r$--compatible and such
that the $U_\alpha$'s cover $M$. A \emph{smooth manifold} is a topological
manifold with a smooth atlas. We use the term \emph{manifold} generally; we specify whether it
is a finite-- or infinite--dimensional manifold only when it is not clear from
context.
If a covering by charts takes their values in a
Banach space $E$, then $E$ is called the \emph{model space} and we say that $M$
is a $C^r$--\emph{Banach manifold}. We remark that one can form a manifold
modeled on any linear space in which one has theory of differential calculus;
we use Banach manifolds so that we can utilize the inverse function theorem. 

Suppose that $f:M\rar N$ where $M, N$ are $C^k$--manifolds. We say $f$ is of
\emph{class}
$C^r$ with $0\leq r\leq k$, and we write $f\in C^r(M, N)$, if for each $\pt\in M$ and a chart $(V,
\psi)$ of $N$ with $f(\pt)\in V$, there is a chart $(U, \vphi)$ of $M$ satisfying
$\pt\in U$, $f(U)\subset V$, and such that the local representation of $f$,
namely $\psi\circ f\circ \vphi^{-1}$, is of class $C^r$.
If $N=\mb{R}$, then $\psi$ can be taken to be the identity map so that 
the local representation is given by $f\circ \vphi^{-1}$.

Each $\pt\in M$ has an associated \emph{tangent space} $T_{\pt} M$, and the
disjoint union of the tangent spaces is the \emph{tangent bundle} $TM =
\coprod_{\pt\in M}T_{\pt} M$.
The \emph{co-tangent space to $M$} at $\pt\in M$, denoted $T^\ast_\pt M$, is the
set of all real-valued linear functionals---or, simply, the dual---on the tangent space $T_\pt M$, and the
disjoint union of the co--tangent spaces is the \emph{co--tangent bundle}
$T^\ast M=\coprod_{\pt\in M} T_\pt^\ast M$. 
Both $TM$ and $T^\ast M$ are naturally smooth manifolds~\cite[Thm.~3.3.10 and
Ch.~5.2 resp.]{abraham:1988aa}. 

For a vector space $E$ we define the vector space of continuous
$(r+s)$--multilinear maps $T_s^r(E)=L^{r+s}(E^\ast, \ldots, E^\ast, E,
\ldots, E; \mb{R})$ with $s$ copies of $E$ and $r$ copes of $E^\ast$ and where $E^\ast$ denotes the dual. We say
elements of $T_s^r(E)$ are \emph{tensors} on $E$, \emph{contravariant} of order
$r$ and \emph{covariant} of order $s$. 
Further, we use the notation $T^r_s(M)$ to
denote the \emph{vector bundle of tensors contravariant of order $r$ and
covariant of order $s$}~\cite[Def.~5.2.9]{abraham:1988aa}. In this notation,
$T^1_0(M)$ is identified with the tangent bundle $TM$ and $T^0_1(M)$ with the
cotangent bundle $T^\ast M$.

 Suppose $f:M\rar N$ is a mapping of one
 manifold into another, and $\pt\in M$, then by means of charts we can interpret
 the derivative of $f$ on each chart at $\pt$ as a linear mapping 
   $df(\pt):T_\pt M\rar T_{f(\pt)}N.$
 When $N=\mb{R}$, the collection of such
maps defines a \emph{$1$--form} $df:M\rar T^\ast M$. More generally, a
$1$--form is a continuous map $\omega:M\rar T^\ast M$ satisfying $\pi\circ
\omega=\text{Id}_M$ where $\pi:T^\ast M\rar M$ is the natural projection mapping
$\omega(p)\in T^\ast_p M$ to $p\in M$.

A point $\pt \in M$ is said to be a \emph{critical point} of a map $f\in
C^r(M,\mb{R})$, $r\geq 2$ if $df(\pt)=0$. 
At a critical point $\pt\in M$, 
there is a uniquely determined continuous, symmetric, bilinear form (termed the
\emph{Hessian})
$d^2f(\pt)\in T_2^0(M)$ such that $d^2f(\pt)$
is defined for all $v,w\in T_uM$ by
  $d^2(f\circ \vphi^{-1})(\vphi(\pt))(v_\vphi, w_\vphi)$
where $\vphi$ is any product chart at $u$ and $v_\vphi, w_\vphi$ are the local representations
of $v,w$ respectively~\cite[Prop. in \S 7]{palais:1963aa}. We say $d^2f(\pt)$ is
\emph{positive semi--definite} if there exists $\alpha\geq 0$ such that for any chart
$\vphi$,
\begin{equation}
  d^2(f\circ \vphi^{-1})(\vphi(u))(v,v)\geq \alpha \|v\|^2, \ \ \forall \ v\in
  T_{\vphi(u)}E.
  \label{eq:posdef}
\end{equation}
If $\alpha>0$, then we say $d^2f(\pt)$ is \emph{positive--definite}.
Both $\omega(\pt)=0$ and positive definiteness are invariant
with respect to the choice of coordinate chart. 

Given a Banach space $E$ and a bounded, symmetric bilinear form $B$ on $E$, we
say that $B$ is \emph{non--degenerate} if the linear map $A:E\rar E^\ast$
defined by $A(v)(w)=B(v,w)$ is a linear isomorphism of $E$ onto $E^\ast$,
otherwise $B$ is \emph{degenerate}.
A critical point $\pt$ of $f$ is called \emph{non--degenerate} if the Hessian of $f$ at
$\pt$ is non--degenerate~\cite[Def. in \S 7]{palais:1963aa}. Degeneracy is independent of the
choice of coordinate chart.

Consider smooth manifolds $M_1, \ldots, M_n$. 
The product space $\prod_{i=1}^n M_i = M_1\times \cdots\times M_n$ is naturally a smooth
manifold~\cite[Def.~3.2.4]{abraham:1988aa}. 
In particular, there is an atlas on $M_1\times \cdots\times M_n$ composed of
\emph{product charts} $(U_1\times \cdots\times U_n, \vphi_1\times\cdots\times
\vphi_n)$ where $(U_i, \vphi_i)$ is a chart on $M_i$ for $i\in \{1,\ldots, n\}$.
We use the notation $\bigtimes_{i=1}^n\vphi_i=\vphi_1\times \cdots \times
\vphi_n$ and $\prod_{i=1}^nU_i=U_1\times \cdots \times U_n$. 

 There is a canonical isomorphism at each point such that the cotangent bundle
 of the product manifold splits:
  \begin{equation}
    T^\ast_{(\pt_1, \ldots,\pt_n)}(M_1\times\cdots\times M_n)\cong
    T^\ast_{\pt_1} M_1\oplus\cdots\oplus
    T^\ast_{\pt_n} M_n
    \label{eq:canon}
  \end{equation}
  where $\oplus$ denotes the direct sum of vector spaces.
 There are natural bundle maps 
 \begin{equation}
 \psi_{M_i}:T^\ast(M_1\times\cdots\times M_n)\rar
 T^\ast(M_1\times\cdots\times M_n)  
   \label{eq:bundlemaps}
 \end{equation}
  annihilating the all the components other
 than those corresponding to $M_i$ of an element in
 the 
 cotangent bundle for each $i\in\{1, \ldots, n\}$. In particular, 
 $\psi_{M_i}(\omega_1, \ldots, \omega_n)=(0, \ldots, 0, \omega_i, 0, \ldots,
 0)$
 where $\omega=(\omega_1, \ldots, \omega_n)\in T^\ast_u(M_1\times\cdots\times M_n)$ and $0$ is
 the zero functional in $T^\ast_{u_j} M_j$ for each $j\neq i$.
 
Let $M=M_1\times \cdots \times M_n$. 
 Given a point $u=(u_1, \ldots, u_n)\in M$, then $\iota_{u}^j: M_j\rar M$ is the natural
 inclusion map where $\iota_{u}^j(\mu)=(u_1, \ldots, u_{j-1}, \mu, u_{j+1},
 \ldots, u_n)$. Suppose we have a function
 $f:M\rar\mb{R}$. Then the derivatives $D_if(u)$ of the map $\mu_i\mapsto f(u_1, \ldots,
 u_{i-1}, \mu_i, u_{i+1}, \ldots, u_n)$ where $\mu_i\in M_i$ for each $i\in\{1, \ldots, n\}$ 
 are called the
 \emph{partial derivatives} of $f$ at $u\in
 M$~\cite[Prop.~2.4.12]{abraham:1988aa}. They are 
 given by $D_if(u)(v_i)=df(u)(\bar{v}_i)$
 where $v_i\in T_{u_i}M_i$ and $\bar{v}_i=(0, \ldots, 0, v_i, 0, \ldots, 0)\in T_uM$.
 Indeed, $d\iota_u^i:T_{u_i}M\rar T_u M$ is a map such that
 $d\iota_u^i(u_i)(v_i)=\bar{v}_i$. Hence, by the chain rule, we have
 $D_if(u)=d(f\circ \iota_{u}^i)(u_i)=df(u)\circ d\iota_u^i.$
 Further, we have that for $v=(v_1, \ldots, v_n)$,
$df(u)(v)=\sum_{i=1}^n D_if(u)(v_i).$ 
For second--order partial derivatives, we use the notation $D_{ij}^2f(\pt)=D_i(D_j f)(\pt)$. 

%

\ifCLASSOPTIONcaptionsoff
  \newpage
\fi

\bibliographystyle{IEEEtran}
\bibliography{2014IEEETACv1}

\end{document}